\documentclass[9pt,shortpaper,twoside,web]{IEEEtran}
\usepackage{cite}

\ifCLASSINFOpdf
  \usepackage[pdftex]{graphicx}
  \usepackage{epstopdf}
  \usepackage{subfigure}
  \graphicspath{{./}{./figures/}}
\else
\fi
\usepackage{xcolor}

\usepackage{amsmath}
\usepackage{amssymb}
\allowdisplaybreaks
\newcommand*{\qed}{\hfill\ensuremath{\square}}%
\newcommand*\diff{\mathop{}\!\mathrm{d}}%

\newtheorem{theorem}{Theorem}[section]
\newtheorem{corollary}{Corollary}[theorem]
\newtheorem{lemma}[theorem]{Lemma}

\newtheorem{assum}[theorem]{Assumption}
\newtheorem{remark}[theorem]{Remark}
\newtheorem{definition}[theorem]{Definition}

\begin{document}
%
\title{Input-to-State Stability of a Clamped-Free Damped String in the Presence of Distributed and Boundary Disturbances}
%
%
%

\author{Hugo~Lhachemi, David~Saussi{\'e}, Guchuan~Zhu, Robert~Shorten
\thanks{Hugo Lhachemi is with the School of Electrical and Electronic Engineering, University College Dublin, Dublin, Ireland (e-mail: hugo.lhachemi@ucd.ie).
 
David Saussi{\'e} and Guchuan Zhu are with the Electrical Engineering Department, Polytechnique Montreal, Montreal, Canada (email: \{d.saussie,guchuan.zhu\}@polymtl.ca).

Robert Shorten is with the School of Electrical and Electronic Engineering, University College Dublin, Dublin, Ireland, and also with the Dyson School of Design Engineering, Imperial College London, London, U.K. (e-mail: robert.shorten@ucd.ie).

This publication has emanated from research supported in part by a research grant from Science Foundation Ireland (SFI) under grant number 16/RC/3872 and is co-funded under the European Regional Development Fund and by I-Form industry partners.}
}


%
%

\markboth{Manuscript submitted to IEEE Transaction on Automatic Control}%
{Lhachemi \MakeLowercase{\textit{et al.}}: Input-to-State Stability of a Clamped-Free Damped String in the Presence of Distributed and Boundary Disturbances}
%



\maketitle

\begin{abstract}
This note establishes the Exponential Input-to-State Stability (EISS) property for a clamped-free damped string with respect to distributed and boundary disturbances. While efficient methods for establishing ISS properties for distributed parameter systems with respect to distributed disturbances have been developed during the last decades, establishing ISS properties with respect to boundary disturbances remains challenging. One of the well-known methods for well-posedness analysis of systems with boundary inputs is the use of a lifting operator for transferring the boundary disturbance to a distributed one. However, the resulting distributed disturbance involves time derivatives of the boundary perturbation. Thus, the subsequent ISS estimate depends on its amplitude, and may not be expressed in the strict form of ISS properties. To solve this problem, we show for a clamped-free damped string equation that the projection of the original system trajectories in an adequate Riesz basis can be used to establish the desired EISS property.
\end{abstract}


%
\IEEEpeerreviewmaketitle

\section{Introduction}\label{sec: Introduction}
%
%
%
%

Originally introduced by Sontag for finite dimensional systems~\cite{sontag1989smooth}, Input-to-State Stability (ISS) is one of the central notions in the modern theory of robust control. Specifically, ISS aims at ensuring that disturbances can only induce, in the worst case, a proportional perturbation of the magnitude of the system trajectory. While this notion has been widely studied for finite dimensional systems, its extension to Partial Differential Equations (PDEs), and more generally to infinite dimensional systems, remains challenging~\cite{jacob2018infinite,mironchenko2016restatements,mironchenko2017characterizations,karafyllis2019preview}. 

For systems described by PDEs, there exist essentially two types of perturbations. The first type includes distributed (or in-domain) perturbations, i.e., perturbations acting over the domain. The second type concerns boundary perturbations, i.e., perturbations acting on the boundary of the domain. This second type of perturbation naturally appears in numerous boundary control problems such as heat equations~\cite{Curtain2012}, transport equations~\cite{karafyllis2016iss}, diffusion or diffusive equations~\cite{argomedo2012d}, and vibration of structures~\cite{Curtain2012} with practical applications, e.g., in robotics~\cite{endo2017boundary,henikl2016infinite}, aerospace engineering~\cite{bialy2016adaptive,lhachemi2018boundary,lhachemi2018boundaryAutomatica}, and additive manufacturing~\cite{diagne2015state,heigel2015thermo}.

In the recent literature, many results have been reported regarding the ISS property with respect to distributed disturbances~\cite{mironchenko2016local,dashkovskiy2011local,argomedo2013strict,dashkovskiy2013input,mironchenko2015construction,mazenc2011strict,prieur2012iss}. In contrast, the literature dealing with the establishment of ISS properties with respect to boundary disturbances is less developed~\cite{argomedo2012d,karafyllis2016iss,karafyllis2016input,karafyllis2017iss,mironchenko2017monotonicity}. The main difficulty relies in the fact that boundary disturbances are generally transferred into distributed disturbances by means of a lifting operator in the framework of the boundary control systems~\cite{Curtain2012}. By doing so, the original system with boundary perturbations is made equivalent to a system with exclusively distributed perturbations, for which efficient tools for analyzing the ISS properties exist. However, the resulting distributed perturbation, and consequently the subsequent ISS estimate, involve time derivatives of the boundary perturbation~\cite{Curtain2012}. Thus, it is paramount to obtain an ISS property compliant with the original definition of ISS, which is exclusively expressed in terms of the amplitude of the disturbances.

A possible approach for establishing ISS properties of PDEs with respect to boundary disturbances consists in resorting to an adequate Lyapunov function~\cite{argomedo2012d,tanwani2017disturbance,zheng2017input,zheng2017giorgi,zheng2018iss}. While very efficient, such an approach relies on the practical capability to construct an adequate Lyapunov function, which is generally challenging and highly case-dependent. Alternative approaches relying on functional analysis tools were investigated in~\cite{karafyllis2016iss,jacob2018infinite,jacob2018continuity}. In~\cite{jacob2018infinite,jacob2018continuity}, the ISS property for disturbances evaluated in the uniform norm is obtained for a class of analytic semigroups. The problem is embedded into the extrapolation space while invoking admissible conditions for returning to the original state-space. A different approach that avoids the incursion into the extrapolation space was developed in~\cite{karafyllis2016iss} for the analysis of 1-D parabolic equations. One of the key ideas was to take advantage of the intrinsic properties of the underlying disturbance-free operator. Indeed, as its opposite belongs to the class of Sturm-Liouville operator, it is self-adjoint, and an adequate selection of a sequence of its eigenvectors provides a Hilbert basis of the underlying Hilbert space. Then, by projecting the system trajectories onto this Hilbert basis and taking advantage of the self-adjoint nature of the disturbance-free operator, it was shown that the analysis of the system trajectories reduces to the study of a countably infinite number of Ordinary Differential Equations (ODEs). Each of these ODEs describes the time domain evolution of one coefficient of the system trajectory in the aforementioned Hilbert basis. The ISS property was finally obtained by solving these ODEs and by resorting to Parseval's identity.

The first motivation of this note is to establish the Exponential Input-to-State (EISS) property of a clamped-free damped string in the presence of both distributed and boundary disturbances. Also known as wave equation, the underlying second-order hyperbolic linear PDE occurs in many fields such as mechanics, acoustics, and fluid dynamics. For this reason, its study in the disturbance free case has attracted a lot of attention~\cite{liu2006exponential,xu2008spectrum,guo2010spectral,guo2012spectrum}. The assessment of the finite asymptotic gain of a clamped damped string in both spatial $L^2$ and sup norms has been reported very recently in~\cite{karafyllis2018boundary} for boundary disturbances of class $\mathcal{C}^4$. The EISS property in spatial sup norm of a similar clamped configuration was also obtained in~\cite{karafyllis2018small} via stability analysis of an equivalent hyperbolic–parabolic PDE loop by means of a small-gain approach. The result presented in this note differs from~\cite{karafyllis2018boundary} as we study the EISS property of a clamped-free configuration (that corresponds to a different set of boundary conditions) in the state-space norm, for disturbance signals evaluated in both uniform and $L^2$ norms, and under the weaker regularity assumptions that the boundary disturbances are of class $\mathcal{C}^2$. 

The second motivation is to show that the approach relying on functional analysis tools employed in~\cite{karafyllis2016iss} can be extended to the problem considered in this work. Specifically, the desired EISS property of the clamped-free damped string is obtained through three main steps. First, the well-posedness of the distributed parameter system is assessed in the framework of boundary control systems. Second, the properties of the underlying disturbance-free operator are studied. Unlike the problem studied in~\cite{karafyllis2016iss}, the disturbance-free operator is not self-adjoint and its eigenvectors do not form a Hilbert basis of the underlying Hilbert space. However, it is a Riesz-spectral operator, implying in particular that its eigenvectors form a Riesz basis~\cite{christensen2016introduction}, which is an important generalization of the concept of Hilbert basis. In particular, even if the Parseval's identity does not hold for Riesz bases, a connexion still exists between the norm of a vector and its coefficients in the Riesz basis. Thus, taking finally advantage of the projection of system trajectories over this Riesz basis, as well as the connection between the eigenstructures of a Riesz-spectral operator and its adjoint operator, we show that the analysis of system trajectories reduces to the study of a countably infinite number of ODEs. By doing so, the EISS property can be derived directly from the original system, which allows avoiding the occurrence of the time derivative of the boundary perturbation. 

The remainder of this note is organized as follows. Notations and definitions are introduced in Section~\ref{sec: notations}. The considered clamped-free damped string model and its well-posedness analysis are presented in Section~\ref{sec: problem}. The detailed study of the properties of the underlying disturbance-free operator is completed in Section~\ref{sec: detailed study A0}. Section~\ref{sec: proof main result} is devoted to the establishment of the EISS property of the system in presence of both distributed and boundary perturbations. Finally, some concluding remarks are provided in Section~\ref{sec: conclusion}.

\section{Notations and definitions}\label{sec: notations}
The sets of non-negative integers, integers, real, non-negative real, positive real, and complex numbers are denoted by $\mathbb{N}$, $\mathbb{Z}$, $\mathbb{R}$, $\mathbb{R}_+$, $\mathbb{R}_+^*$, and $\mathbb{C}$, respectively. For any $z \in \mathbb{C}$, $\operatorname{Re}z$ and $\operatorname{Im}z$ denote the real part and the imaginary part of $z$, respectively. For any integer $k\in\mathbb{Z}$, we define $\tilde{k} \triangleq k + 1/2$. We define for $N \in \mathbb{N}$ the following sets: 
\begin{align*}
\mathcal{I}_N & = \{ (k,\epsilon) \,:\, 0 \leq k \leq N ,\, \epsilon \in \{-1,+1\} \} , \\
\mathcal{I}_\infty & = \mathbb{N} \times \{-1,+1\} .
\end{align*} 

For an interval $I\subset\mathbb{R}$ and a normed space $(E,\Vert\cdot\Vert_E)$, $\mathcal{C}^n(I;E)$ (simply $\mathcal{C}^n(I)$ when $E=\mathbb{C}$ endowed with the absolute value) denotes the set of functions $f : I \rightarrow E$ that are $n$ times continuously differentiable. For any $a<b$, we endowed $\mathcal{C}^0([a,b];E)$ with the usual norm $\Vert\cdot\Vert_{\mathcal{C}^0([a,b];E)}$ defined for any $f\in\mathcal{C}^0([a,b];E)$ by 
\begin{equation*}
\Vert f \Vert_{\mathcal{C}^0([a,b];E)} = \underset{t \in [a,b]}{\sup} \Vert f(t) \Vert_E .
\end{equation*} 

The set of square-integrable functions (w.r.t. the Lebesgue measure) over an interval $(a,b) \subset \mathbb{R}$ is denoted by $L^2(a,b)$ and is endowed with its natural inner product $\left< f , g \right>_{L^2(a,b)} = \int_{a}^{b} f(\xi) \overline{g(\xi)} \diff\xi$, providing a structure of Hilbert space. Denoting by $f'$, when it exists, the weak derivative of $f \in L^2(a,b)$, we consider the Sobolev space $H^1(a,b) \triangleq \{ f \in L^2(a,b) \; : \; f' \in L^2(a,b) \}$. Finally, $H_L^1(a,b) \triangleq \{ f \in H^1(a,b) \; : \; f(a)=0 \}$ is endowed with the inner product $\left< f , g \right>_{H_L^1(a,b)} \triangleq \left< f' , g' \right>_{L^2(a,b)}$, providing a structure of Hilbert space.

For a given linear operator $L$, $\mathrm{R}(L)$, $\mathrm{ker}(L)$, and $\rho(L)$ denote its range, its kernel, and its resolvent set, respectively. The set of linear bounded operators $L: E \rightarrow E$ is denoted by $\mathcal{L}(E)$. The time derivative of a complex-valued differentiable function $f: I \rightarrow\mathbb{C}$ is denoted by $\dot{f}$. Denoting by $(\mathcal{H},\left<\cdot,\cdot\right>_\mathcal{H})$ a $\mathbb{C}$-Hilbert space, the time derivative of a $\mathcal{H}$-valued differentiable function $f: I \rightarrow\mathcal{H}$ is denoted by $\mathrm{d} f / \mathrm{d}t$.

Finally, we introduce the following classical definitions.
\begin{definition}[Riesz basis~\cite{christensen2016introduction}]\label{def: Riesz basis}
A sequence $\Phi = \{ \varphi_k ,\; k\in\mathbb{N}\}$ of vectors of $\mathcal{H}$ is a Riesz basis if 1) $\Phi$ is maximal: $\overline{\mathrm{span}_\mathbb{C}(\Phi)} = \mathcal{H}$, i.e., the closure of the vector space spanned by $\Phi$ coincides with the whole space $\mathcal{H}$; 2) there exist $m_R,M_R\in\mathbb{R}_+^*$ such that for any $N\in\mathbb{N}$ and any $a_k \in\mathbb{C}$, 
\begin{equation}\label{eq: Riesz basis inequality}
m_R \sum\limits_{0 \leq k \leq N} |a_k|^2
\leq
\left\Vert \sum\limits_{0 \leq k \leq N} a_{k} \varphi_{k} \right\Vert_\mathcal{H}^2
\leq
M_R \sum\limits_{0 \leq k \leq N} |a_k|^2 . 
\end{equation}
\end{definition}

\begin{definition}[Riesz spectral operator~\cite{Curtain2012}]\label{def: Riesz-spectral operator}
Let $A : D(A) \subset \mathcal{H} \rightarrow \mathcal{H}$ be a linear and closed operator with simple eigenvalues $\lambda_n$ and corresponding eigenvectors $\varphi_n \in D(A)$, $n \in \mathbb{N}$. Operator $A$ is a Riesz-spectral operator if 1) $\left\{ \varphi_n , \; n \in \mathbb{N} \right\}$ is a Riesz basis; 2) the closure of $\{ \lambda_n , \; n \in \mathbb{N} \}$ is totally disconnected, i.e., for any two distinct $a,b \in \overline{ \{ \lambda_n , \; n \in \mathbb{N} \} }$, $[a,b] \not\subset \overline{ \{ \lambda_n , \; n \in \mathbb{N} \} }$.
\end{definition}

\section{Problem Description and Main Result}\label{sec: problem}

\subsection{Problem setting}
Consider a string with Kelvin-Voigt damping~\cite{liu2006exponential,xu2008spectrum,guo2010spectral} and clamped-free boundary conditions described by:
\begin{subequations}
\begin{align}
\dfrac{\partial^2 y}{\partial t^2} - \dfrac{\partial}{\partial x} \left( \alpha \dfrac{\partial y}{\partial x} + \beta \dfrac{\partial^2 y}{\partial t \partial x} \right) & = u , & & \; \mathrm{in}\;\mathbb{R}_{+}\times(0,1) \label{eq: EDP} \\
y(t,0) & = 0 , & & \; t \in \mathbb{R}_+ \label{eq: left BC} \\
\left( \alpha \dfrac{\partial y}{\partial x} + \beta \dfrac{\partial^2 y}{\partial t \partial x} \right)(t,1) & = d(t) , & & \; t \in \mathbb{R}_+ \label{eq: right BC} \\
y(0,x) & = y_0(x) , & & \; x \in (0,1)  \label{eq: CI 1} \\
\dfrac{\partial y}{\partial t} (0,x) & = y_{t0}(x) , & & \; x \in (0,1)  \label{eq: CI 2}
\end{align}
\end{subequations} 
where $\alpha,\beta\in\mathbb{R}_+^*$ are constant parameters. Functions $u\in\mathcal{C}^1(\mathbb{R}_+;L^2(0,1))$ and $d\in\mathcal{C}^2(\mathbb{R}_+)$ represent distributed and boundary disturbances, respectively. Functions $y_0 \in H_L^1(0,1)$ and $y_{t0} \in L^2(0,1)$ are the initial conditions. 

Throughout the paper, we assume that the following assumption holds. Its introduction is motivated by the properties of the underlying operators, as it will be shown in the subsequent developments.

\begin{assum}\label{assumption}
The coefficients $\alpha,\beta\in\mathbb{R}_+^*$ in the system (\ref{eq: EDP}-\ref{eq: CI 2}) are such that
\begin{equation*}
\dfrac{2 \sqrt{\alpha}}{\pi\beta}-\dfrac{1}{2}\notin\mathbb{N} .
\end{equation*}
\end{assum}

To study (\ref{eq: EDP}-\ref{eq: CI 2}), we introduce the functional space:
\begin{equation*}
\mathcal{H} = H_L^1(0,1) \times L^2(0,1) ,
\end{equation*}
which is a Hilbert space when endowed with the inner product defined for all $(x_1,x_2),(\hat{x}_1,\hat{x}_2)\in\mathcal{H}$ by
\begin{equation*}
\left<(x_1,x_2) , (\hat{x}_1,\hat{x}_2) \right>_\mathcal{H} = \int_0^1 \alpha x_1'(\xi) \overline{\hat{x}_1'(\xi)} + x_2(\xi) \overline{\hat{x}_2(\xi)} \diff\xi .
\end{equation*}
Let the operator $\mathcal{A} : D(\mathcal{A}) \rightarrow \mathcal{H}$ be defined by 
\begin{equation*}
\mathcal{A}(x_1,x_2) = (x_2,(\alpha x_1' + \beta x_2')') 
\end{equation*}
over the domain
\begin{align*}
D(\mathcal{A}) = \{ (x_1,x_2) \in \mathcal{H} \; : & \; x_2 \in H_L^1(0,1) , \\
& (\alpha x_1' + \beta x_2' ) \in H^1(0,1) \} .
\end{align*}
Let $\mathcal{B}: D(\mathcal{B}) \rightarrow \mathbb{C}$ be the boundary operator defined by 
\begin{equation*}
\mathcal{B}(x_1,x_2) = (\alpha x_1' + \beta x_2')(1)
\end{equation*}
with $D(\mathcal{B}) = D(\mathcal{A})$. Finally, introducing $U=(0,u)\in\mathcal{C}^1(\mathbb{R}_+;\mathcal{H})$, (\ref{eq: EDP}-\ref{eq: CI 2}) can be written under the following abstract form~\cite{Curtain2012}:
\begin{equation}\label{eq: abstract form}
\left\{\begin{split}
\dfrac{\mathrm{d} X}{\mathrm{d} t}(t) & = \mathcal{A} X(t) + U(t) & ,\;  t \geq 0\\
\mathcal{B}X(t) & = d(t) & ,\; t \geq 0 \\
X(0) & = X_0 
\end{split}\right.
\end{equation}
with the state $X(t)=(y(t,\cdot),y_t(t,\cdot))$ and the initial condition $X_0=(y_0,y_{t0})$.

\subsection{Well-posedness}\label{subsec: well-posedness}
We introduce the disturbance-free operator $\mathcal{A}_0$ defined over the domain $D(\mathcal{A}_0) \triangleq D(\mathcal{A}) \cap \mathrm{ker}(\mathcal{B})$ by $\mathcal{A}_0 \triangleq \left. \mathcal{A} \right|_{D(\mathcal{A}_0)}$. Straightforward computations show that $\mathcal{A}_0$ generates a $C_0$-semigroup of contractions. Indeed, $\mathcal{A}_0$ is dissipative since a simple integration by parts yields for any $(x_1,x_2)\in D(\mathcal{A}_0)$,
\begin{equation*}
\operatorname{Re} \left( \left< \mathcal{A}_0(x_1,x_2) , (x_1,x_2) \right>_\mathcal{H} \right) = - \beta \int_0^1 |x_2'(\xi)|^2 \diff\xi \leq 0 .
\end{equation*}
Furthermore, direct computations show that $\mathcal{A}_0$ is invertible and is defined for any $(x_1,x_2)\in\mathcal{H}$ by
\begin{equation*}
\mathcal{A}_0^{-1}(x_1,x_2) = \left( - \dfrac{\beta}{\alpha} x_1 - \dfrac{1}{\alpha} \int_0^{(\cdot)} \int_{\xi_1}^1 x_2(\xi_2) \diff\xi_2 \diff\xi_1 , x_1\right).
\end{equation*}
Finally, Poincar{\'e} and Cauchy-Schwarz inequalities~\cite{Hardy1952} imply that $\mathcal{A}_0^{-1} \in \mathcal{L}(\mathcal{H})$. The application of the Lumer-Phillips theorem~\cite{Pazy2012,luo2012stability} yields the desired result, i.e., $\mathcal{A}_0$ generates a $C_0$-semigroup of contractions $T$.

In order to conclude on the well-posedness of the abstract system (\ref{eq: abstract form}) for an initial condition $X_0 \in D(\mathcal{A})$ and perturbations $u\in\mathcal{C}^1(\mathbb{R}_+;L^2(0,1))$ and $d\in\mathcal{C}^2(\mathbb{R}_+)$ such that $\mathcal{B}X_0 = d(0)$, it is sufficient to check that the abstract system satisfies the definition of a \emph{Boundary control system}~\cite[Def. 3.3.2]{Curtain2012}. Introducing the lifting operator $B : \mathbb{C} \rightarrow \mathcal{H}$ defined for any $d\in\mathbb{C}$ by $Bd = (f_d,0)$ with $f_d(x) \triangleq (d/\alpha)x$ for all $x \in [0,1]$, we have $\mathrm{R}(B)\subset D(\mathcal{A})$, $\mathcal{A}B=0_{\mathcal{L}(\mathbb{C},\mathcal{H})}$, and $\mathcal{B}B = I_\mathbb{C}$. Thus the abstract system (\ref{eq: abstract form}) is well-posed for any $X_0 \in D(\mathcal{A})$, $u\in\mathcal{C}^1(\mathbb{R}_+;L^2(0,1))$, and $d\in\mathcal{C}^2(\mathbb{R}_+)$ such that $\mathcal{B}X_0 = d(0)$~\cite[Th 3.1.3 ; Th. 3.3.3]{Curtain2012}. Furthermore, $X \in \mathcal{C}^0(\mathbb{R}_+;D(\mathcal{A}))\cap \mathcal{C}^1(\mathbb{R}_+;\mathcal{H})$ is the classical solution\footnote{The function $X$ is a classical solution of (\ref{eq: abstract form}) if $X \in \mathcal{C}^0(\mathbb{R}_+;D(\mathcal{A}))\cap \mathcal{C}^1(\mathbb{R}_+;\mathcal{H})$ and satisfies (\ref{eq: abstract form}) for all $t \geq 0$.} of (\ref{eq: abstract form}) if and only if $V = X - Bd \in \mathcal{C}^0(\mathbb{R}_+;D(\mathcal{A}_0))\cap \mathcal{C}^1(\mathbb{R}_+;\mathcal{H})$ is the classical solution of the following abstract system:
\begin{equation*}
\left\{\begin{split}
\dfrac{\diff V}{\diff t}(t) & = \mathcal{A}_0 V(t) - B \dot{d}(t) + U(t) , & \; t \geq 0  \\
V(0) & = V_0 
\end{split}\right.
\end{equation*}
where $V_0 = X_0 - Bd(0) \in D(\mathcal{A}_0)$. Thus, by direct integration, the solution of (\ref{eq: abstract form}) is given for $t \geq 0$ by
\begin{align}
X(t) 
& = T(t) \left( X_0 - B d(0) \right) + B d(t) \label{eq: explicit solution} \\
& \phantom{=}\, + \int_0^t T(t-\tau) \left\{ -B \dot{d}(\tau) + U(\tau) \right\} \diff\tau . \nonumber
\end{align}

\begin{remark}
It is pointed out that a weak version of the ISS property, implying the norm of time derivative $\dot{d}$ of the boundary disturbance $d$, can be easily obtained from (\ref{eq: explicit solution}). However, the traditional definition of ISS is stronger because it is only limited to the amplitude of the boundary disturbance, and not its time derivatives. The objective of this paper is to establish such an ISS estimate for (\ref{eq: EDP}-\ref{eq: CI 2}), only with respect to the magnitude of the perturbations $d$ and $u$.
\end{remark}

\subsection{Main result}
Throughout the paper, let $k_0\in\mathbb{N}$ be defined by
\begin{equation}\label{eq: def k0}
k_0 \triangleq \left\lceil \dfrac{2\sqrt{\alpha}}{\pi\beta} - \dfrac{1}{2} \right\rceil \geq 0 ,
\end{equation}
where $\lceil \cdot \rceil$ denotes the ceiling function.

The main result of the paper regarding the ISS property of the trajectories of the abstract system (\ref{eq: abstract form}) is stated in the following theorem. 

\begin{theorem}\label{th: main result}
For any initial condition $X_0 \in D(\mathcal{A})$ and any disturbances $u\in\mathcal{C}^1(\mathbb{R}_+;L^2(0,1))$ and $d\in\mathcal{C}^2(\mathbb{R}_+)$ such that $\mathcal{B}X_0 = d(0)$, the abstract system (\ref{eq: abstract form}) has a unique classical solution $X \in \mathcal{C}^0(\mathbb{R}_+;D(\mathcal{A}))\cap \mathcal{C}^1(\mathbb{R}_+;\mathcal{H})$. Furthermore, under Assumption~\ref{assumption}, the system is EISS with respect to disturbances in both uniform and $L^2$ norms in the sense that there exist constants $C_0,C_1,C_2,C_3,C_4\in\mathbb{R}_+^*$, independent of $X_0$, $u$, and $d$, such that for all $t \geq 0$,
\begin{align}
\Vert X(t) \Vert_\mathcal{H} \leq & C_0 e^{- \kappa_0 t} \Vert X_0 \Vert_\mathcal{H} + C_1 \Vert d \Vert_{\mathcal{C}^0([0,t])} \label{eq: main result - ISS} \\ 
& + C_2 \Vert u \Vert_{\mathcal{C}^0([0,t];L^2(0,1))} , \nonumber
\end{align}
and
\begin{align}
\Vert X(t) \Vert_\mathcal{H} \leq & C_0 e^{- \kappa_0 t} \Vert X_0 \Vert_\mathcal{H} + C_3 \Vert d \Vert_{L^2(0,t)} \label{eq: main result - ISS 2} \\ 
& + C_4 \Vert u \Vert_{L^2((0,t)\times(0,1))} , \nonumber
\end{align}
where
\begin{equation}\label{eq: main result - growth bound}
\kappa_0 = 
\left\{\begin{split}
\min \left( \dfrac{\beta\pi^2}{8} , \dfrac{\alpha}{\beta} \right) 
& \;\; \mathrm{if} \; k_0 \geq 1 ; \\
\dfrac{\alpha}{\beta} 
& \;\; \mathrm{if} \; k_0 = 0 ,
\end{split}\right.
\end{equation}
and is such that $\omega_0 = - \kappa_0 < 0$ is the growth bound of $T(t)$.
\end{theorem}

The proof of Theorem~\ref{th: main result} is presented in Section~\ref{sec: proof main result} after studying the spectral properties of $\mathcal{A}_0$ in Section~\ref{sec: detailed study A0}.

\begin{remark}
As $\omega_0 = - \kappa_0$ is the growth bound of $T$, the convergence rate of the exponential term in (\ref{eq: main result - ISS}-\ref{eq: main result - ISS 2}) is tight in the sense that $\kappa_0$ cannot be replaced in (\ref{eq: main result - ISS}) or (\ref{eq: main result - ISS 2}) by any $\kappa > \kappa_0$ such that the ISS estimate still holds true.
\end{remark}

\section{Study of the properties of $\mathcal{A}_0$}\label{sec: detailed study A0}
The objective of this section is to demonstrate that $\mathcal{A}_0$ is a Riesz-spectral operator while characterizing the eigenvalues and eigenfunctions of both $\mathcal{A}_0$ and $\mathcal{A}_0^*$.

\subsection{Characterization of the spectral properties of $\mathcal{A}_0$}

\begin{lemma}\label{lemma: A0 eigenvalues and eigenvectors}
The eigenvalues of $\mathcal{A}_0$ are simple and are given by $\{\lambda_{k,\epsilon} , \; k\in\mathbb{N} , \; \epsilon\in\{-1,+1\}\}$ where
\begin{equation}\label{eq: A0 eigenvalues}
\lambda_{k,\epsilon} = \left\{\begin{split}
-\dfrac{\tilde{k}^2\beta\pi^2}{2} + \epsilon i  \dfrac{\tilde{k}\pi\sqrt{4\alpha - \tilde{k}^2\beta^2\pi^2}}{2}
& , & 0 \leq k \leq k_0 - 1 ; \\ 
-\dfrac{\tilde{k}^2\beta\pi^2}{2} + \epsilon \dfrac{\tilde{k}\pi\sqrt{\tilde{k}^2\beta^2\pi^2 - 4\alpha}}{2}
& , & k \geq k_0 ,
\end{split}\right.
\end{equation}
with $\tilde{k} \triangleq k + 1/2$. Furthermore, the associated eigenspaces are given by $\mathrm{ker}(\mathcal{A}_0 - \lambda_{k,\epsilon} I_\mathcal{H}) = \mathrm{span}_\mathbb{C}(\phi_{k,\epsilon})$ with
\begin{equation}\label{eq: A0 eigenvectors}
\phi_{k,\epsilon} = \dfrac{1}{\lambda_{k,\epsilon}} \left( \sin(\tilde{k}\pi\cdot) , \lambda_{k,\epsilon} \sin(\tilde{k}\pi\cdot) \right) .
\end{equation}
\end{lemma} 

\textbf{Proof.} Let $\lambda\in\mathbb{C}$ and $(x_1,x_2) \in D(\mathcal{A}_0)\backslash\{0\}$ be such that $\mathcal{A}_0 (x_1,x_2) = \lambda (x_1,x_2)$, i.e., $x_1(0)=x_2(0)=(\alpha x_1' + \beta x_2')(1)=0$ with $x_2 = \lambda x_1 $ and $(\alpha x_1' + \beta x_2')' = \lambda x_2$. As $\lambda \neq -\alpha/\beta$ (otherwise we would have $x_1=x_2=0$), we deduce that 
\begin{equation*}
x_2 = \lambda x_1 , \qquad x_1'' = \dfrac{\lambda^2}{(\alpha + \lambda\beta)} x_1 ,
\end{equation*}
with $x_1'(1)=0$. Denoting by $r(\lambda)\in\mathbb{C}$ one of the two distinct\footnote{Because $\lambda^2/(\alpha + \lambda\beta)=0$ implies $\lambda =0$ which yields $x_1=x_2=0$.} square-roots of $\lambda^2/(\alpha + \lambda\beta)$, there exist  $a,b\in\mathbb{C}$ such that 
\begin{equation*}
x_1(\xi) = a e^{r(\lambda)\xi} + b e^{-r(\lambda)\xi} .
\end{equation*}
From the boundary conditions we get $b=-a$ and $a r(\lambda) (e^{r(\lambda)} + e^{-r(\lambda)}) = 0$. As $\lambda \neq 0$ and because we are looking for non-trivial solutions (i.e., such that $(x_1,x_2) \neq 0$), $a r(\lambda) \neq 0$ whence $e^{2r(\lambda)}=-1$. We obtain that $2 r(\lambda) \equiv i \pi \; (2i\pi)$, i.e., there exists $k\in\mathbb{Z}$ such that $r(\lambda) = i\tilde{k}\pi$. From the definition of $r(\lambda)$, we deduce that
\begin{equation*}
r(\lambda)^2 = \dfrac{\lambda^2}{(\alpha + \lambda\beta)} = - \tilde{k}^2 \pi^2.
\end{equation*}
Thus $P_k(\lambda) = 0$ where $P_k = X^2 + \tilde{k}^2 \beta \pi^2 X + \tilde{k}^2 \alpha \pi^2 \in \mathbb{R}[X]$. As $P_{-k-1} = P_{k}$, the study for $k\in\mathbb{Z}$ reduces to $k\in\mathbb{N}$. The discriminant of $P_k$ is given by $\mathrm{disc}(P_k) = \tilde{k}^2 \pi^2 ( \tilde{k}^2 \beta^2 \pi^2 - 4 \alpha )$. Based on Assumption~\ref{assumption}, $\mathrm{disc}(P_k) \neq 0$, thus $\mathrm{disc}(P_k) > 0$ for $k \geq k_0$ and $\mathrm{disc}(P_k) < 0$ for $0 \leq k \leq k_0 -1$, providing the eigenvalues $\lambda_{k,\epsilon}$ given by (\ref{eq: A0 eigenvalues}). Finally the associated eigenvectors are characterized by
\begin{equation*}
x_1(\xi) = a ( e^{i\tilde{k}\pi\xi} - e^{-i\tilde{k}\pi\xi}) = 2 a i \sin(\tilde{k}\pi\xi) ,
\end{equation*}
and $x_2 = \lambda_{k,\epsilon} x_1$, providing (\ref{eq: A0 eigenvectors}). \qed

In order to work with unitary eigenvectors, we introduce $\Phi_{k,\epsilon} \triangleq \phi_{k,\epsilon} / \Vert \phi_{k,\epsilon} \Vert_\mathcal{H}$, where a straightforward integration shows that 
\begin{equation}\label{eq: norm phi,k,epsilon}
\Vert \phi_{k,\epsilon} \Vert_\mathcal{H} = \dfrac{1}{\sqrt{2}} \sqrt{ 1 + \dfrac{\tilde{k}^2\alpha\pi^2}{\vert \lambda_{k,\epsilon} \vert^2}} .
\end{equation}
Finally, we denote $\Phi = \{ \Phi_{k,\epsilon} , \; k\in\mathbb{N} , \; \epsilon \in \{-1,+1\} \}$.

For the upcoming developments, we establish certain equalities, inequalities, and asymptotic behaviours for the eigenvalues of $\mathcal{A}_0$. First, note that 
\begin{equation}\label{eq: prod two ev}
\forall k \geq 0 , \;\; \lambda_{k,-1} \lambda_{k,+1} = \tilde{k}^2\alpha\pi^2,
\end{equation}
\begin{equation}\label{eq: maj real complex}
\forall 0 \leq k \leq k_0 - 1, \; \forall \epsilon\in\{-1,+1\}, \;\; \operatorname{Re}\lambda_{k,\epsilon} \leq - \dfrac{\beta \pi^2}{8} .
\end{equation}
Furthermore, as $\sqrt{1+x} \leq 1+x/2$ for all $x \geq -1$, 
\begin{equation}\label{eq: maj real real}
\forall k \geq k_0 , \;\; \lambda_{k,-1} \leq \lambda_{k,+1} \leq - \dfrac{\alpha}{\beta} .
\end{equation}
To study the asymptotic behaviours, we consider $k \geq k_0$, giving
\begin{equation}\label{eq: aymptotic behaviour lambda,k,-1}
\lambda_{k,-1} \underset{k \rightarrow +\infty}{\sim} - k^2 \beta \pi^2 ,
\end{equation}
and
\begin{align}
\lambda_{k,+1} 
& = \dfrac{\tilde{k}^2 \beta \pi^2}{2} \left( -1 + \sqrt{ 1 - \dfrac{4\alpha}{\tilde{k}^2 \beta^2 \pi^2}} \right) \nonumber \\
& = \dfrac{\tilde{k}^2 \beta \pi^2}{2} \left( -1 + \left\{ 1 - \dfrac{2\alpha}{\tilde{k}^2\beta^2 \pi^2} + o(k^{-2}) \right\} \right) \nonumber \\
& = -\alpha/\beta + o(1) \nonumber \\
& \underset{k \rightarrow +\infty}{\xrightarrow{\hspace*{0.75cm}}} -\alpha/\beta . \label{eq: aymptotic behaviour lambda,k,+1}
\end{align}

\subsection{Characterization and properties of $\mathcal{A}_0^*$}
In the subsequent developments, we will use the adjoint operator $\mathcal{A}_0^*$ and in particular the connections between the eigenstructures of $\mathcal{A}_0$ and $\mathcal{A}_0^*$. This is motivated by the fact that $\Phi$ is not a Hilbert basis for $\mathcal{H}$ (more details provided latter in Subsection~\ref{subsec: A0 riesz operator}).

\begin{lemma}
The adjoint operator $\mathcal{A}_0^*$ is defined over the domain
\begin{align*}
D(\mathcal{A}_0^*) = \{ (x_1,x_2) \in \mathcal{H} \; : & \; x_2 \in H_L^1(0,1) , \\
& (\alpha x_1' - \beta x_2' ) \in H^1(0,1) , \\
& (\alpha x_1' - \beta x_2' )(1) = 0 \} ,
\end{align*}
by
\begin{equation*}
\mathcal{A}_0^* (x_1,x_2) =
\left(
-x_2 , - (\alpha x_1' - \beta x_2' )' 
\right) .
\end{equation*}
\end{lemma}

\textbf{Proof.} As $\mathcal{A}_0^{-1} \in \mathcal{L}(\mathcal{H})$, $(\mathcal{A}_0^{-1})^* \in \mathcal{L}(\mathcal{H})$ and $(\mathcal{A}_0^*)^{-1} = (\mathcal{A}_0^{-1})^*$~\cite[Th. III.5.30]{kato2013perturbation}. Integration by parts and application of Fubini theorem yields for any $(x_1,x_2) \in \mathcal{H}$, 
\begin{equation*}
(\mathcal{A}_0^{-1})^{*}(x_1,x_2) = \left( -\dfrac{\beta}{\alpha} x_1 + \dfrac{1}{\alpha} \int_0^{(\cdot)} \int_{\xi_1}^1 x_2(\xi_2) \diff\xi_2 \diff\xi_1 , - x_1 \right) .
\end{equation*}
The inversion of $(\mathcal{A}_0^{-1})^{*}$ gives the claimed result. \qed

\begin{lemma}\label{lemma: A0* eigenvalues and eigenvectors}
The eigenvalues of $\mathcal{A}_0^*$ are given by $\{\mu_{k,\epsilon} , \; k\in\mathbb{N} , \; \epsilon\in\{-1,+1\}\}$ where $\mu_{k,\epsilon} = \overline{\lambda_{k,\epsilon}}$. Furthermore, the associated eigenspaces are given by $\mathrm{ker}(\mathcal{A}_0^* - \mu_{k,\epsilon} I_\mathcal{H}) = \mathrm{span}_\mathbb{C}(\psi_{k,\epsilon})$ with
\begin{equation}\label{eq: A0* eigenvectors}
\psi_{k,\epsilon} = \dfrac{1}{\mu_{k,\epsilon}} \left( - \sin(\tilde{k}\pi\cdot) , \mu_{k,\epsilon} \sin(\tilde{k}\pi\cdot) \right) .
\end{equation}
\end{lemma} 

\textbf{Proof.} Let $\mu\in\mathbb{C}$ and $(x_1,x_2) \in D(\mathcal{A}_0^*)\backslash\{0\}$ be such that $\mathcal{A}_0^* (x_1,x_2) = \mu (x_1,x_2)$, i.e., $x_1(0)=x_2(0)=(\alpha x_1' - \beta x_2')(1)=0$ with $- x_2 = \mu x_1 $ and $-(\alpha x_1' - \beta x_2')' = \mu x_2$. We deduce that 
\begin{equation*}
x_2 = - \mu x_1 , \qquad x_1'' = \dfrac{\mu^2}{(\alpha + \mu\beta)} x_1 ,
\end{equation*}
with $x_1'(1)=0$. Therefore $x_1$ satisfies the same differential equation as the one in the proof of Lemma~\ref{lemma: A0 eigenvalues and eigenvectors} where $x_2 = \lambda x_1$ is replaced by $x_2 = - \mu x_1$. Thus the claimed conclusion follows from the proof of Lemma~\ref{lemma: A0 eigenvalues and eigenvectors}. \qed

For any $(k_1,\epsilon_1) \neq (k_2,\epsilon_2)$, 
\begin{align*}
\lambda_{k_1,\epsilon_1} \left< \Phi_{k_1,\epsilon_1} , \psi_{k_2,\epsilon_2} \right>_\mathcal{H}
& =
\left< \mathcal{A}_0 \Phi_{k_1,\epsilon_1} , \psi_{k_2,\epsilon_2} \right>_\mathcal{H} \\
& =
\left< \Phi_{k_1,\epsilon_1} , \mathcal{A}_0^* \psi_{k_2,\epsilon_2} \right>_\mathcal{H} \\& =
\lambda_{k_2,\epsilon_2} \left< \Phi_{k_1,\epsilon_1} , \psi_{k_2,\epsilon_2} \right>_\mathcal{H} . 
\end{align*}
Thus, as\footnote{It directly follows from Assumption~\ref{assumption} and the fact that $\lambda_{k,\epsilon}^2/(\alpha+\lambda_{k,\epsilon}\beta) = - \tilde{k}^2\pi^2$.} $\lambda_{k_1,\epsilon_1} \neq \lambda_{k_2,\epsilon_2} $, $\left< \Phi_{k_1,\epsilon_1} , \psi_{k_2,\epsilon_2} \right>_\mathcal{H} = 0$. Furthermore, based on $\lambda_{k,\epsilon} \overline{\mu_{k,\epsilon}} = \lambda_{k,\epsilon}^2$ and (\ref{eq: prod two ev}), a direct integration yields
\begin{equation*}
\left< \Phi_{k,\epsilon} , \psi_{k,\epsilon} \right>_\mathcal{H}
= \dfrac{1}{2 \Vert\phi_{k,\epsilon}\Vert_\mathcal{H}} \left( 1 - \dfrac{\lambda_{k,-\epsilon}}{\lambda_{k,\epsilon}} \right) \neq 0 ,
\end{equation*}
because Assumption~\ref{assumption} implies $\lambda_{k,\epsilon} \neq \lambda_{k,-\epsilon}$.
Thus, introducing 
\begin{equation*}
\Psi_{k,\epsilon} \triangleq \dfrac{1}{\;\overline{\left< \Phi_{k,\epsilon} , \psi_{k,\epsilon} \right>_\mathcal{H}}\;} \psi_{k,\epsilon} ,
\end{equation*}
and letting $\Psi = \{ \Psi_{k,\epsilon} , \; k\in\mathbb{N} , \; \epsilon \in \{-1,+1\} \}$, the set of eigenvectors $\Phi$ of $\mathcal{A}_0$ is biorthogonal to the set of eigenvectors $\Psi$ of $\mathcal{A}_0^*$ in the sense that $\left< \Phi_{k_1,\epsilon_1} , \Psi_{k_2,\epsilon_2} \right>_\mathcal{H} = \delta_{(k_1,\epsilon_1),(k_2,\epsilon_2)}$.

\subsection{$\mathcal{A}_0$ is a Riesz-Spectral Operator}\label{subsec: A0 riesz operator}
We show that $\Phi$ is a Riesz basis of $\mathcal{H}$ (see Definition~\ref{def: Riesz basis}) and $\mathcal{A}_0$ is a Riesz-spectral operator (see Definition~\ref{def: Riesz-spectral operator}).

\subsubsection{$\Phi$ is maximal}
Let us first introduce the following technical lemma whose proof is provided in Appendix.

\begin{lemma}\label{lemma: maximal families}
Both $\{ \cos(\tilde{k}\pi\cdot) , \; k\in\mathbb{N}\}$ and $\{ \sin(\tilde{k}\pi\cdot) , \; k\in\mathbb{N}\}$ are maximal in $L^2(0,1)$.
\end{lemma}

Then, the following result holds true.

\begin{lemma}\label{lemma: Phi maximal}
$\Phi$ is maximal in $\mathcal{H}$.
\end{lemma}

\textbf{Proof.} 
Let $z = (z_1,z_2) \in \mathcal{H}$ be such that $\left< \Phi_{k,\epsilon} , z \right>_\mathcal{H} = 0$ for all $k \in \mathbb{N}$ and $\epsilon \in \{-1,+1\}$. Then,
\begin{equation*}
\alpha \tilde{k} \pi \left< \cos(\tilde{k}\pi\cdot) , z_1' \right>_{L^2(0,1)} + \lambda_{k,\epsilon} \left< \sin(\tilde{k}\pi\cdot) , z_2 \right>_{L^2(0,1)} = 0 ,
\end{equation*}
from which we obtain that, for all $k \in \mathbb{N}$,
\begin{equation*}
\begin{bmatrix}
\alpha \tilde{k} \pi & \lambda_{k,-1} \\
\alpha \tilde{k} \pi & \lambda_{k,+1}
\end{bmatrix}
\begin{bmatrix}
\left< \cos(\tilde{k}\pi\cdot) , z_1' \right>_{L^2(0,1)} \\
\left< \sin(\tilde{k}\pi\cdot) , z_2 \right>_{L^2(0,1)}
\end{bmatrix}
= 0 .
\end{equation*}
Based on Assumption~\ref{assumption}, $\alpha \tilde{k} \pi ( \lambda_{k,+1} - \lambda_{k,-1} ) \neq 0$, which implies the invertibility of the $2 \times 2$ matrix. Therefore, we have 
\begin{equation*}
\forall k \in \mathbb{N} , \qquad
\left< \cos(\tilde{k}\pi\cdot) , z_1' \right>_{L^2(0,1)} = 
\left< \sin(\tilde{k}\pi\cdot) , z_2 \right>_{L^2(0,1)} = 0 .
\end{equation*}
Hence, Lemma~\ref{lemma: maximal families} ensures that $z_1'=z_2=0$. As $z_1(0)=0$, we conclude that $z=0$. \qed

\subsubsection{$\Phi$ is a Riesz basis}
Direct integrations show that for any non negative integers $k_1 \neq k_2$ and $\epsilon_1,\epsilon_2 \in \{-1,+1\}$, 
\begin{equation}\label{eq: inner product Phi different k}
\left< \Phi_{k_1,\epsilon_1} , \Phi_{k_2,\epsilon_2} \right>_\mathcal{H} = 0.
\end{equation}
Nevertheless, $\Phi$ is not a Hilbert basis because for any $k\in\mathbb{N}$ and $\epsilon \in \{-1,+1\}$, 
\begin{equation*}
\left< \Phi_{k,\epsilon} , \Phi_{k,-\epsilon} \right>_\mathcal{H} = \dfrac{1}{2 \Vert \phi_{k,\epsilon} \Vert_\mathcal{H} \Vert \phi_{k,-\epsilon} \Vert_\mathcal{H}} \left( 1 + \dfrac{\tilde{k}^2 \alpha \pi^2}{\lambda_{k,\epsilon}\overline{\lambda_{k,-\epsilon}}} \right) \neq 0 .
\end{equation*}
However, we have the following result.

\begin{lemma}\label{lem: Riesz basis}
$\Phi$ is a Riesz basis.
\end{lemma}

\textbf{Proof.} 
Based on Lemma~\ref{lemma: Phi maximal} it is sufficient to show (\ref{eq: Riesz basis inequality}). For any $N \in \mathbb{N}$ and any $a_{k,\epsilon} \in \mathbb{C}$, we infer from (\ref{eq: inner product Phi different k}) that
\begin{align}
& \left\Vert \sum\limits_{(k,\epsilon) \in \mathcal{I}_N} a_{k,\epsilon} \Phi_{k,\epsilon} \right\Vert_\mathcal{H}^2 \nonumber \\
& \quad = \sum\limits_{(k_1,\epsilon_1) \in \mathcal{I}_N} \sum\limits_{(k_2,\epsilon_2) \in \mathcal{I}_N}  a_{k_1,\epsilon_1} \overline{a_{k_2,\epsilon_2}} \left< \Phi_{k_1,\epsilon_1} , \Phi_{k_2,\epsilon_2} \right>_\mathcal{H} \nonumber \\
& \quad = \sum\limits_{0 \leq k \leq N} S_k , \label{eq: norm evaluation for riesz basis}
\end{align}
where 
\begin{equation*}
S_k = |a_{k,-1}|^2 + |a_{k,+1}|^2 + 2 \operatorname{Re} \left( a_{k,-1} \overline{a_{k,+1}} \left< \Phi_{k,-1} , \Phi_{k,+1} \right>_\mathcal{H} \right) .
\end{equation*}
We evaluate the term $\left< \Phi_{k,-1} , \Phi_{k,+1} \right>_\mathcal{H}$ as follows:
\begin{equation*}
\left< \Phi_{k,-1} , \Phi_{k,+1} \right>_\mathcal{H} = \dfrac{ 1 + \dfrac{\tilde{k}^2 \alpha \pi^2}{\lambda_{k,-1}\overline{\lambda_{k,+1}}} }{ \sqrt{ 1 + \dfrac{\tilde{k}^2 \alpha \pi^2}{|\lambda_{k,-1}|^2} } \sqrt{ 1 + \dfrac{\tilde{k}^2 \alpha \pi^2}{|\lambda_{k,+1}|^2} } }.
\end{equation*}

We first consider the case $k \geq k_0$. As $\lambda_{k,-1}\overline{\lambda_{k,+1}} = \lambda_{k,-1}\lambda_{k,+1}=\tilde{k}^2\alpha\pi^2$, we have that
\begin{equation*}
\left< \Phi_{k,-1} , \Phi_{k,+1} \right>_\mathcal{H} = \dfrac{ 2 }{ \sqrt{ 2 + \dfrac{|\lambda_{k,-1}|^2 + |\lambda_{k,+1}|^2}{\tilde{k}^2 \alpha \pi^2} } }.
\end{equation*}
Based on (\ref{eq: A0 eigenvalues}), 
\begin{equation*}
\dfrac{|\lambda_{k,-1}|^2 + |\lambda_{k,+1}|^2}{\tilde{k}^2 \alpha \pi^2} 
= \dfrac{\tilde{k}^2 \beta^2 \pi^2}{\alpha} - 2
\geq \dfrac{\tilde{k}_0^2 \beta^2 \pi^2}{\alpha} - 2 ,
\end{equation*}
yielding for all $k \geq k_0$, 
\begin{equation*}
\left| \left< \Phi_{k,-1} , \Phi_{k,+1} \right>_\mathcal{H} \right|
\leq \dfrac{4\sqrt{\alpha}}{(2 k_0+1) \beta \pi} < 1 ,
\end{equation*}
where the last inequality holds true because, based on the definition (\ref{eq: def k0}) of $k_0$ and Assumption~\ref{assumption}, $k_0 > 2 \sqrt{\alpha} / (\beta\pi) -1/2$.

We now consider the case $0 \leq k \leq k_0 -1$ when $k_0 \geq 1$. A similar computation shows that
\begin{equation*}
\left\vert \left< \Phi_{k,-1} , \Phi_{k,+1} \right>_\mathcal{H}  \right\vert
= \dfrac{\tilde{k} \beta \pi}{2 \sqrt{\alpha}}
\leq \dfrac{(2 k_0 - 1)\beta \pi}{4 \sqrt{\alpha}} 
< 1 ,
\end{equation*}
where the last inequality holds true because, based on the definition (\ref{eq: def k0}) of $k_0$ and Assumption~\ref{assumption}, $k_0 < 2 \sqrt{\alpha} / (\beta\pi) +1/2$.

Thus, introducing 
\begin{equation*}
C \triangleq 
\left\{\begin{split}
\max \left( \dfrac{4 \sqrt{\alpha}}{(2 k_0 +1) \beta \pi} , \dfrac{(2 k_0 - 1)\beta \pi}{4 \sqrt{\alpha}} \right) 
& \;\; \mathrm{if} \; k_0 \geq 1 ; \\
\dfrac{4 \sqrt{\alpha}}{\beta \pi} 
& \;\; \mathrm{if} \; k_0 = 0 ,
\end{split}\right.
\end{equation*}
we obtain that $C \in (0,1)$ and $\left|\left< \Phi_{k,-1} , \Phi_{k,+1} \right>_\mathcal{H}\right| \leq C$ for all $k \geq 0$. This yields 
\begin{equation*}
\left\vert \operatorname{Re} \left( a_{k,-1} \overline{a_{k,+1}} \left< \Phi_{k,-1} , \Phi_{k,+1} \right>_\mathcal{H} \right) \right\vert
\leq C \vert a_{k,-1} \vert \vert a_{k,+1} \vert.
\end{equation*}
Consequently, we have
\begin{align*}
S_k
& \leq \vert a_{k,-1} \vert^2 + \vert a_{k,+1} \vert^2 + 2 C \vert a_{k,-1} \vert \vert a_{k,+1} \vert \\
& \leq (1-C) \left( \vert a_{k,-1} \vert^2 + \vert a_{k,+1} \vert^2 \right) + C \left( \vert a_{k,-1} \vert + \vert a_{k,+1} \vert \right)^2 \\
& \leq (1+C) \left( \vert a_{k,-1} \vert^2 + \vert a_{k,+1} \vert^2 \right) ,
\end{align*}
and
\begin{align*}
S_k
& \geq \vert a_{k,-1} \vert^2 + \vert a_{k,+1} \vert^2 - 2 C \vert a_{k,-1} \vert \vert a_{k,+1} \vert \\
& \geq (1-C) \left( \vert a_{k,-1} \vert^2 + \vert a_{k,+1} \vert^2 \right) + C \left( \vert a_{k,-1} \vert - \vert a_{k,+1} \vert \right)^2 \\
& \geq (1-C) \left( \vert a_{k,-1} \vert^2 + \vert a_{k,+1} \vert^2 \right) .
\end{align*}
Combining the two inequalities above with (\ref{eq: norm evaluation for riesz basis}), we obtain the desired result:
\begin{align}
m_R \sum\limits_{(k,\epsilon) \in \mathcal{I}_N} |a_{k,\epsilon}|^2
& \leq
\left\Vert \sum\limits_{(k,\epsilon) \in \mathcal{I}_N} a_{k,\epsilon} \Phi_{k,\epsilon} \right\Vert_\mathcal{H}^2 \nonumber \\
& \leq
M_R \sum\limits_{(k,\epsilon) \in \mathcal{I}_N} |a_{k,\epsilon}|^2 \label{eq: Riesz basis inequality for campled-free damped string}
\end{align}
with $m_R = 1 - C > 0$ and $M_R = 1 + C > 0$. As $m_R$ and $M_R$ are constants independent of $N \in \mathbb{N}$ and $a_k \in \mathbb{C}$, the claimed conclusion holds true. \qed

\begin{remark}
The constants $m_R = 1 - C$ and $M_R = 1 + C$ provide a tight version of (\ref{eq: Riesz basis inequality for campled-free damped string}). Indeed, it follows from the proof that there exists $k \in \{ k_0 - 1 , k_0 \} $ such that $\left|\left< \Phi_{k,-1} , \Phi_{k,+1} \right>_\mathcal{H}\right| = C$. Thus $\left< \Phi_{k,-1} , \Phi_{k,+1} \right>_\mathcal{H} = C e^{i\theta}$ for some $\theta\in[0,2\pi)$. Considering $a_{k,-1} = \overline{a_{k,+1}} = e^{-i\theta/2}$ we obtain $S_k = 2(1+C) = (1+C) \left( \vert a_{k,-1} \vert^2 + \vert a_{k,+1} \vert^2 \right)$. Conversely, with $a_{k,-1} = - \overline{a_{k,+1}} = e^{-i\theta/2}$ we obtain $S_k = 2(1-C) = (1-C) \left( \vert a_{k,-1} \vert^2 + \vert a_{k,+1} \vert^2 \right)$.
\end{remark}

As $\Phi$ is a Riesz basis biorthogonal to $\Psi$, we obtain from the general theory on Riesz basis~\cite{christensen2016introduction} that for all $z\in\mathcal{H}$, 
\begin{equation}\label{eq: Riesz basis decomposition equality}
z = \sum\limits_{(k,\epsilon) \in \mathcal{I}_\infty} \left<z,\Psi_{k,\epsilon}\right>_\mathcal{H} \Phi_{k,\epsilon} ,
\end{equation}
and
\begin{align}
(1-C) & \sum\limits_{(k,\epsilon) \in \mathcal{I}_\infty} |\left<z,\Psi_{k,\epsilon}\right>_\mathcal{H}|^2
\leq
\left\Vert z \right\Vert_\mathcal{H}^2 \nonumber \\
& \leq
(1+C) \sum\limits_{(k,\epsilon) \in \mathcal{I}_\infty} |\left<z,\Psi_{k,\epsilon}\right>_\mathcal{H}|^2 . \label{eq: Riesz basis decomposition inequality}
\end{align}

\subsubsection{$\mathcal{A}_0$ is a Riesz-Spectral Operator}
We can now state the main result of this section.

\begin{lemma}
The operator $\mathcal{A}_0$ is a Riesz-spectral operator generating an exponentially stable $C_0$-semigroup with growth $\omega_0=-\kappa_0 < 0$ where $\kappa_0$ is given by (\ref{eq: main result - growth bound}) 
\end{lemma}

\textbf{Proof.} 
We directly deduce from the fact that $\mathcal{A}_0$ generates a $C_0$-semigroup and from Lemmas~\ref{lemma: A0 eigenvalues and eigenvectors} and~\ref{lem: Riesz basis} that $\mathcal{A}_0$ is a Riesz-spectral operator. Thus, its growth bound $\omega_0$ satisfies~\cite[Th. 2.3.5]{Curtain2012}: 
\begin{equation*}
\omega_0 = \underset{(k,\epsilon)\in\mathcal{I}_\infty}{\sup} \operatorname{Re} \lambda_{k,\epsilon} .
\end{equation*}
Based on (\ref{eq: maj real complex}-\ref{eq: maj real real}), $\operatorname{Re} \lambda_{k,\epsilon} \leq - \kappa_0$ where $\kappa_0$ is given by (\ref{eq: main result - growth bound}). If $k_0 \geq 1$, (\ref{eq: maj real complex}) becomes an equality for $k=0$. Furthermore, as (\ref{eq: aymptotic behaviour lambda,k,+1}) holds, this yields $\omega_0 = - \kappa_0$. \qed

\section{Proof of the EISS property}\label{sec: proof main result}

We can now prove the main result of this note.  

\textbf{Proof of Theorem~\ref{th: main result}.} Let $X_0 \in D(\mathcal{A})$, $u\in\mathcal{C}^1(\mathbb{R}_+;L^2(0,1))$, and $d\in\mathcal{C}^2(\mathbb{R}_+)$ such that $\mathcal{B}X_0 = d(0)$. Let $X=(x_1,x_2) \in \mathcal{C}^0(\mathbb{R}_+;D(\mathcal{A}))\cap \mathcal{C}^1(\mathbb{R}_+;\mathcal{H})$ be the unique classical solution of the abstract system (\ref{eq: abstract form}). Based on (\ref{eq: Riesz basis decomposition equality}-\ref{eq: Riesz basis decomposition inequality}),
\begin{equation}\label{eq: main result initiale maj}
\forall t\geq 0, \;\;
\left\Vert X(t) \right\Vert_\mathcal{H}^2 
\leq
(1+C) \sum\limits_{(k,\epsilon)\in\mathcal{I}_\infty} | c_{k,\epsilon}(t) |^2 ,
\end{equation}
where $c_{k,\epsilon} \triangleq \left<X,\Psi_{k,\epsilon}\right>_\mathcal{H} \in \mathcal{C}^1(\mathbb{R}_+)$. With $\Psi_{k,\epsilon} \triangleq (\Psi_{k,\epsilon}^1,\Psi_{k,\epsilon}^2)$, we have for all $t \geq 0$, 
\begin{align}
\dot{c}_{k,\epsilon}(t)
= & \left< \dfrac{\mathrm{d} X}{\mathrm{d} t}(t),\Psi_{k,\epsilon}\right>_\mathcal{H} \nonumber \\
= & \left< \mathcal{A}X(t)+U(t),\Psi_{k,\epsilon}\right>_\mathcal{H} \nonumber \\
= & \left< (x_2(t),(\alpha x_1' + \beta x_2')'(t) + u(t)),(\Psi_{k,\epsilon}^1,\Psi_{k,\epsilon}^2)\right>_\mathcal{H} \nonumber \\
= & \int_0^1 \alpha x_2'(t) \overline{\Psi_{k,\epsilon}^{1'}} + \{(\alpha x_1' + \beta x_2')'(t) + u(t)\} \overline{\Psi_{k,\epsilon}^2} \diff\xi \nonumber \\
= & \int_0^1 \alpha x_2'(t) \overline{\Psi_{k,\epsilon}^{1'}} \diff\xi + \left[ (\alpha x_1' +\beta x_2')(t) \overline{\Psi_{k,\epsilon}^2}
 \right]_{\xi=0}^{\xi=1} \nonumber \\
 & - \int_0^1 (\alpha x_1' +\beta x_2')(t) \overline{\Psi_{k,\epsilon}^{2'}} \diff\xi + \int_0^1 u(t) \overline{\Psi_{k,\epsilon}^{2}} \diff\xi \nonumber \\
= & \int_0^1 \alpha x_1'(t) \overline{(-\Psi_{k,\epsilon}^{2})'} \diff\xi + \int_0^1 x_2'(t) \overline{ \{ \alpha \Psi_{k,\epsilon}^{1'} - \beta \Psi_{k,\epsilon}^{2'} \} } \diff\xi \nonumber \\
& + d(t) \overline{\Psi_{k,\epsilon}^2(1)} + \int_0^1 u(t) \overline{\Psi_{k,\epsilon}^{2}} \diff\xi \nonumber \\
= & \int_0^1 \alpha x_1'(t) \overline{(-\Psi_{k,\epsilon}^{2})'} \diff\xi + \left[ x_2(t) \overline{ \{ \alpha \Psi_{k,\epsilon}^{1'} - \beta \Psi_{k,\epsilon}^{2'} \} } \right]_{\xi=0}^{\xi=1} \nonumber \\
& - \int_0^1 x_2(t) \overline{ \{ \alpha \Psi_{k,\epsilon}^{1'} - \beta \Psi_{k,\epsilon}^{2'} \}' } \diff\xi + d(t) \overline{\Psi_{k,\epsilon}^2(1)} \nonumber \\ 
& + \int_0^1 u(t) \overline{\Psi_{k,\epsilon}^{2}} \diff\xi \nonumber \\
= & \left< X(t), \mathcal{A}_0^* \Psi_{k,\epsilon}\right>_\mathcal{H} + d(t) \overline{\Psi_{k,\epsilon}^2(1)} + \int_0^1 u(t) \overline{\Psi_{k,\epsilon}^{2}} \diff\xi . \label{eq: proof main result - key idea}
\end{align}
As $\mathcal{A}_0^* \Psi_{k,\epsilon} = \mu_{k,\epsilon} \Psi_{k,\epsilon} = \overline{\lambda_{k,\epsilon}} \Psi_{k,\epsilon}$, this yields for any $t \geq 0$,
\begin{equation*}
\dot{c}_{k,\epsilon}(t) = \lambda_{k,\epsilon} c_{k,\epsilon}(t) + d(t) \overline{\Psi_{k,\epsilon}^2(1)} + \int_0^1 u(t) \overline{\Psi_{k,\epsilon}^{2}} \diff\xi ,
\end{equation*}
which gives after integration:
\begin{align}
c_{k,\epsilon}(t) 
& = e^{\lambda_{k,\epsilon} t} c_{k,\epsilon}(0) 
+ \int_0^t e^{\lambda_{k,\epsilon} (t-\tau)} d(\tau) \overline{\Psi_{k,\epsilon}^2(1)} \diff\tau \nonumber \\
& \phantom{=} + \int_0^t e^{\lambda_{k,\epsilon} (t-\tau)} \int_0^1 u(\tau) \overline{\Psi_{k,\epsilon}^{2}} \diff\xi \diff\tau . \label{eq: expression of c,k,eps}
\end{align}
We estimate the three terms on the right-hand side of (\ref{eq: expression of c,k,eps}) as follows. First,
\begin{equation}\label{eq: upper bound 1}
\left| e^{\lambda_{k,\epsilon} t} c_{k,\epsilon}(0) \right|
\leq e^{\operatorname{Re}\lambda_{k,\epsilon} t} \left| c_{k,\epsilon}(0) \right|
\leq e^{- \kappa_0 t} \left| c_{k,\epsilon}(0) \right| .
\end{equation}
Second, introducing $\gamma_{k,\epsilon} \triangleq \vert \Psi_{k,\epsilon}^2(1) / \operatorname{Re}\lambda_{k,\epsilon} \vert$,
\begin{align}
& \left| \int_0^t e^{\lambda_{k,\epsilon} (t-\tau)} d(\tau) \overline{\Psi_{k,\epsilon}^2(1)} \diff\tau \right| \nonumber \\
& \qquad\leq \gamma_{k,\epsilon} \int_0^t - \operatorname{Re}\lambda_{k,\epsilon} e^{\operatorname{Re}\lambda_{k,\epsilon} (t-\tau)} \diff\tau \, \Vert d \Vert_{\mathcal{C}^0([0,t])} \nonumber \\
& \qquad\leq \gamma_{k,\epsilon} \left( 1 - e^{\operatorname{Re}\lambda_{k,\epsilon} t} \right) \Vert d \Vert_{\mathcal{C}^0([0,t])} \nonumber \\
& \qquad\leq \gamma_{k,\epsilon} \Vert d \Vert_{\mathcal{C}^0([0,t])} , \label{eq: upper bound 2}
\end{align}
with
\begin{align*}
\gamma_{k,\epsilon}
& = \left| \dfrac{\psi_{k,\epsilon}^2(1)}{\operatorname{Re}(\lambda_{k,\epsilon}) \left< \Phi_{k,\epsilon} , \psi_{k,\epsilon} \right>_\mathcal{H}} \right| \\
& = \dfrac{2 \Vert \phi_{k,\epsilon} \Vert_\mathcal{H}}{\left\vert \operatorname{Re}(\lambda_{k,\epsilon}) \left( 1 - \dfrac{\lambda_{k,-\epsilon}}{\lambda_{k,\epsilon}} \right) \right\vert} . 
\end{align*}
Finally, by using Cauchy-Schwarz inequality,  
\begin{align*}
& \left\vert \int_0^t e^{\lambda_{k,\epsilon} (t-\tau)} \int_0^1 u(\tau) \overline{\Psi_{k,\epsilon}^{2}} \diff\xi \diff\tau \right\vert \\
\leq & \dfrac{1}{\left| \operatorname{Re} \lambda_{k,\epsilon} \right|} \int_0^t - \operatorname{Re}\lambda_{k,\epsilon} e^{\operatorname{Re}\lambda_{k,\epsilon} (t-\tau)}   \int_0^1 \left\vert u(\tau) \overline{\Psi_{k,\epsilon}^{2}} \right\vert \diff\xi \diff\tau \\
\leq & \dfrac{\Vert \Psi_{k,\epsilon}^{2} \Vert_{L^2(0,1)}}{\left| \operatorname{Re} \lambda_{k,\epsilon} \right|} \int_0^t - \operatorname{Re}\lambda_{k,\epsilon} e^{\operatorname{Re}\lambda_{k,\epsilon} (t-\tau)}  \Vert u(\tau) \Vert_{L^2(0,1)} \diff\tau \\
\leq & \dfrac{\Vert \Psi_{k,\epsilon}^{2} \Vert_{L^2(0,1)}}{\left| \operatorname{Re} \lambda_{k,\epsilon} \right|} \left( 1 - e^{\operatorname{Re}\lambda_{k,\epsilon} t} \right) \Vert u \Vert_{\mathcal{C}^0([0,t];L^2(0,1))} \\
\leq & \dfrac{\Vert \Psi_{k,\epsilon}^{2} \Vert_{L^2(0,1)}}{\left| \operatorname{Re} \lambda_{k,\epsilon} \right|} \Vert u \Vert_{\mathcal{C}^0([0,t];L^2(0,1))} ,
\end{align*}
and as 
\begin{align*}
\Vert \Psi_{k,\epsilon}^{2} \Vert_{L^2(0,1)}
& = \dfrac{1}{\left| \left< \Phi_{k,\epsilon} , \psi_{k,\epsilon} \right>_\mathcal{H} \right|} \sqrt{ \int_0^1 \sin^2(\tilde{k} \pi \xi) \diff\xi } \\
& = \dfrac{\sqrt{2} \Vert \phi_{k,\epsilon} \Vert_\mathcal{H}}{\left\vert  1 - \dfrac{\lambda_{k,-\epsilon}}{\lambda_{k,\epsilon}} \right\vert} \\
& = \gamma_{k,\epsilon} \vert\operatorname{Re}\lambda_{k,\epsilon}\vert / \sqrt{2} ,
\end{align*}
this yields
\begin{align}
& \left\vert \int_0^t e^{\lambda_{k,\epsilon} (t-\tau)} \int_0^1 u(\tau) \overline{\Psi_{k,\epsilon}^{2}} \diff\xi \diff\tau \right\vert \nonumber \\ 
& \qquad\qquad\qquad \leq 
\dfrac{\sqrt{2}}{2} \gamma_{k,\epsilon} \Vert u \Vert_{\mathcal{C}^0([0,t];L^2(0,1))} . \label{eq: upper bound 3}
\end{align}
Putting together (\ref{eq: expression of c,k,eps}) with the inequalities (\ref{eq: upper bound 1}-\ref{eq: upper bound 3}), this yields for all $t \geq 0$,
\begin{align*}
\left| c_{k,\epsilon}(t) \right|
\leq &
e^{-\kappa_0 t} \left| c_{k,\epsilon}(0) \right|
+ \gamma_{k,\epsilon} \Vert d \Vert_{\mathcal{C}^0([0,t])} \\
& + \dfrac{\sqrt{2}}{2} \gamma_{k,\epsilon} \Vert u \Vert_{\mathcal{C}^0([0,t];L^2(0,1))}.
\end{align*}
As $(a+b+c)^2 \leq 3 (a^2 + b^2 + c^2)$ for all $a,b,c\in\mathbb{R}$, 
\begin{align}
\left| c_{k,\epsilon}(t) \right|^2
\leq &
3 e^{-2\kappa_0 t} \left| c_{k,\epsilon}(0) \right|^2
+ 3 \gamma_{k,\epsilon}^2 \Vert d \Vert_{\mathcal{C}^0([0,t])}^2 \nonumber \\
& + \dfrac{3}{2} \gamma_{k,\epsilon}^2 \Vert u \Vert_{\mathcal{C}^0([0,t];L^2(0,1))}^2. \label{eq: main result maj c,k,eps(t)}
\end{align}
We need to check that $\gamma_{k,\epsilon}^2$ is a summable sequence. To do so, considering $k \geq k_0$, $\operatorname{Re} \lambda_{k,\epsilon} = \lambda_{k,\epsilon}$, which gives along with (\ref{eq: A0 eigenvalues}) and (\ref{eq: norm phi,k,epsilon})
\begin{equation*}
\gamma_{k,\epsilon}
= \dfrac{\sqrt{2}}{\tilde{k} \pi \sqrt{\tilde{k}^2 \beta^2 \pi^2 - 4 \alpha}} \sqrt{1 + \dfrac{\tilde{k}^2 \alpha \pi^2}{\vert \lambda_{k,\epsilon} \vert^2}} .
\end{equation*}
Based on (\ref{eq: aymptotic behaviour lambda,k,-1}-\ref{eq: aymptotic behaviour lambda,k,+1}), the following asymptotic behaviours hold
\begin{equation*}
\gamma_{k,+1}
\underset{k \rightarrow +\infty}{\sim} \dfrac{1}{k \pi} \sqrt{\dfrac{2}{\alpha}} ,
\qquad
\gamma_{k,-1}
\underset{k \rightarrow +\infty}{\sim} \dfrac{\sqrt{2}}{k^2 \beta \pi^2} ,
\end{equation*}
assessing that $\gamma_{k,\epsilon}$ is a square summable sequence. Therefore, we can define the constant $\gamma\in\mathbb{R}_+$ by 
\begin{equation*}
\gamma^2 \triangleq \sum\limits_{(k,\epsilon)\in\mathcal{I}_\infty} \gamma_{k,\epsilon}^2 < \infty .
\end{equation*}
Noting that, based on (\ref{eq: Riesz basis decomposition equality}-\ref{eq: Riesz basis decomposition inequality}),
\begin{equation*}
\sum\limits_{(k,\epsilon)\in\mathcal{I}_\infty} | c_{k,\epsilon}(0) |^2
\leq 
\dfrac{1}{1-C} \Vert X_0 \Vert_\mathcal{H}^2 ,
\end{equation*}
we obtain by using (\ref{eq: main result maj c,k,eps(t)}) into (\ref{eq: main result initiale maj}) that
\begin{align*}
\left\Vert X(t) \right\Vert_\mathcal{H}^2
\leq &
3 \dfrac{1+C}{1-C} e^{-2\kappa_0 t} \Vert X_0 \Vert_\mathcal{H}^2
+ 3 (1+C) \gamma^2 \Vert d \Vert_{\mathcal{C}^0([0,t])}^2 \\
& + \dfrac{3}{2} (1+C) \gamma^2 \Vert u \Vert_{\mathcal{C}^0([0,t];L^2(0,1))}^2 .
\end{align*}
To conclude, it is sufficient to note that $\sqrt{a+b} \leq \sqrt{a} + \sqrt{b}$ for all $a,b\in\mathbb{R}_+$, which yields
\begin{align*}
\left\Vert X(t) \right\Vert_\mathcal{H}
\leq &
\sqrt{3 \dfrac{1+C}{1-C}} e^{-\kappa_0 t} \Vert X_0 \Vert_\mathcal{H}
+ \gamma \sqrt{3 (1+C)} \Vert d \Vert_{\mathcal{C}^0([0,t])} \\
& + \gamma \sqrt{\dfrac{3}{2} (1+C)} \Vert u \Vert_{\mathcal{C}^0([0,t];L^2(0,1))} .
\end{align*}
Thus, the claimed ISS estimate (\ref{eq: main result - ISS}) holds with 
\begin{equation*}
C_0 = \sqrt{3 \dfrac{1+C}{1-C}}, \; C_1 = \gamma \sqrt{3 (1+C)} , \; C_2 = \gamma \sqrt{\dfrac{3}{2} (1+C)} .
\end{equation*}

To prove the second ISS estimate (\ref{eq: main result - ISS 2}), we substitute the estimations (\ref{eq: upper bound 2}-\ref{eq: upper bound 3}) with
\begin{equation*}
\left| \int_0^t e^{\lambda_{k,\epsilon} (t-\tau)} d(\tau) \overline{\Psi_{k,\epsilon}^2(1)} \diff\tau \right| 
\leq \sqrt{\dfrac{\vert \operatorname{Re} \lambda_{k,\epsilon} \vert}{2}} \gamma_{k,\epsilon} \Vert d \Vert_{L^2(0,t)} ,
\end{equation*}
and
\begin{align*}
& \left\vert \int_0^t e^{\lambda_{k,\epsilon} (t-\tau)} \int_0^1 u(\tau) \overline{\Psi_{k,\epsilon}^{2}} \diff\xi \diff\tau \right\vert \nonumber \\ 
& \qquad \qquad \leq 
\dfrac{1}{2} \sqrt{\vert \operatorname{Re} \lambda_{k,\epsilon} \vert} \gamma_{k,\epsilon} \Vert u \Vert_{L^2((0,t)\times(0,1))} ,
\end{align*}
where the Cauchy-Schwartz inequality has been used. This yields for all $t \geq 0$,
\begin{align*}
\left| c_{k,\epsilon}(t) \right|
\leq &
e^{-\kappa_0 t} \left| c_{k,\epsilon}(0) \right|
+ \sqrt{\dfrac{\vert \operatorname{Re} \lambda_{k,\epsilon} \vert}{2}} \gamma_{k,\epsilon} \Vert d \Vert_{L^2(0,t)} \\
& + \dfrac{1}{2} \sqrt{\vert \operatorname{Re} \lambda_{k,\epsilon} \vert} \gamma_{k,\epsilon} \Vert u \Vert_{L^2((0,t)\times(0,1))} .
\end{align*}
Noting that, for any $\epsilon \in \{-1,+1\}$,
\begin{equation*}
\sqrt{\vert \operatorname{Re} \lambda_{k,\epsilon} \vert} \gamma_{k,\epsilon}
\underset{k \rightarrow +\infty}{\sim} \dfrac{1}{k \pi} \sqrt{\dfrac{2}{\beta}} ,
\end{equation*}
$\sqrt{\vert \operatorname{Re} \lambda_{k,\epsilon} \vert} \gamma_{k,\epsilon}$ is a square summable sequence and we can define the constant $\gamma'\in\mathbb{R}_+$ by 
\begin{equation*}
\gamma'^2 \triangleq \sum\limits_{(k,\epsilon)\in\mathcal{I}_\infty} \vert \operatorname{Re} \lambda_{k,\epsilon} \vert \gamma_{k,\epsilon}^2 < \infty .
\end{equation*}
Following the same procedure used above to demonstrate the ISS estimate (\ref{eq: main result - ISS}), we obtain for all $t \geq 0$,
\begin{align*}
\left\Vert X(t) \right\Vert_\mathcal{H}
\leq &
\sqrt{3 \dfrac{1+C}{1-C}} e^{-\kappa_0 t} \Vert X_0 \Vert_\mathcal{H}
+ \gamma' \sqrt{\dfrac{3}{2} (1+C)} \Vert d \Vert_{L^2(0,t)} \\
& + \dfrac{\gamma'}{2} \sqrt{3 (1+C)} \Vert u \Vert_{L^2((0,t)\times(0,1))} .
\end{align*}
Thus, introducing the constants $C_3,C_4\in\mathbb{R}_+^*$ defined by 
\begin{equation*}
C_3 = \gamma' \sqrt{\dfrac{3}{2} (1+C)} , \; C_4 = \dfrac{\gamma'}{2} \sqrt{3 (1+C)} ,
\end{equation*}
the second claimed ISS estimate (\ref{eq: main result - ISS 2}) holds.
\qed

\begin{remark}
The key idea in the proof of the main result lies in the computation of (\ref{eq: proof main result - key idea}). In the disturbance free case, i.e., $d = 0$ and $u = 0$, one has $X \in \mathcal{C}^0(\mathbb{R}_+;D(\mathcal{A}_0))\cap \mathcal{C}^1(\mathbb{R}_+;\mathcal{H})$. Then, because $\Psi_{k,\epsilon} \in D(\mathcal{A}_0^*)$, we obtain as a direct consequence of the definition of the adjoint operator that
\begin{align*}
\dot{c}_{k,\epsilon}(t)
= \left< \dfrac{\mathrm{d} X}{\mathrm{d} t}(t),\Psi_{k,\epsilon}\right>_\mathcal{H}
& = \left< \mathcal{A}_0 X(t),\Psi_{k,\epsilon}\right>_\mathcal{H} \\
& = \left< X(t) , \mathcal{A}_0^* \Psi_{k,\epsilon}\right>_\mathcal{H} ,
\end{align*}
which coincides with (\ref{eq: proof main result - key idea}) in the disturbance free case. In the disturbed case, the computation (\ref{eq: proof main result - key idea}) is nothing but the heuristic computation of the adjoint operator while letting appear 1) an extra non zero boundary condition, via the integrations by parts, due to the boundary disturbance $d$ ; 2) an integral term due to the distributed disturbance $U$.
\end{remark}

\begin{remark}
Putting together (\ref{eq: Riesz basis decomposition equality}) and (\ref{eq: expression of c,k,eps}), one can get an explicit formula of the system trajectory $X$ in function of $X_0$, $d$, $U$, and the eigenstructures of operators $\mathcal{A}_0$ and $\mathcal{A}_0^*$. 
\end{remark}

We deduce, as a direct consequence of the ISS estimates (\ref{eq: main result - ISS}-\ref{eq: main result - ISS 2}) and of the semigroup property of (\ref{eq: abstract form}), the following asymptotic behaviour.
\begin{corollary}
Under the notations and assumptions of Theorem~\ref{th: main result}, assume that one of the two following conditions holds:
\begin{itemize}
\item the perturbations are vanishing in the sense that $\vert d(t) \vert \underset{t \rightarrow +\infty}{\longrightarrow} 0$ and $\Vert u(t) \Vert_{L^2(0,1)} \underset{t \rightarrow +\infty}{\longrightarrow} 0$ ;
\item the perturbations are of finite energy, i.e., $d \in L^2(\mathbb{R}_+)$ and $u \in L^2(\mathbb{R}_+;L^2(0,1)) \cong L^2(\mathbb{R}_+\times(0,1))$,
\end{itemize}
then $\Vert X(t) \Vert_\mathcal{H} \underset{t \rightarrow +\infty}{\longrightarrow} 0$.
\end{corollary}

\section{Conclusion}\label{sec: conclusion}
This paper established the property of Exponential Input-to-State Stability (EISS) for a clamped-free damped string with respect to distributed and boundary disturbances. The adopted approach does not rely on the construction of an adequate Lyapunov function but takes advantage of functional analysis tools. Specifically, by projecting the system trajectories onto a Riesz basis of the underlying Hilbert space formed by the eigenvectors of the disturbance-free operator, the EISS property was derived directly on the original system, avoiding the appearance of the time derivative of the boundary perturbation.


%
%
%

\appendix[Proof of Lemma~\ref{lemma: maximal families}]

From the Fourier series theory~\cite{korner1989fourier}, the set of functions $\{e^{i k \pi \cdot} , \; k\in\mathbb{Z}\}$ is a Hilbert basis of $L^2(-1,1)$ endowed with $\left<f,g\right>_{L^2(-1,1)} = \dfrac{1}{2}\int_{-1}^1 f(\xi) \overline{g(\xi)} \diff\xi$. Let $f \in L^2(-1,1)$ and consider $\hat{f} = e^{-i\pi \cdot /2} f \in L^2(-1,1)$. As
\begin{equation*}
\hat{f} = \sum\limits_{k \in \mathbb{Z}} \left< \hat{f} , e^{ik\pi\cdot} \right>_{L^2(-1,1)} e^{ik\pi\cdot} ,
\end{equation*}
and $\vert e^{i\pi \cdot /2} \vert = 1$, then
\begin{equation*}
f = e^{i\pi \cdot /2} \hat{f} = \sum\limits_{k \in \mathbb{Z}} \left< \hat{f} , e^{ik\pi\cdot} \right>_{L^2(-1,1)} e^{i\tilde{k}\pi\cdot} .
\end{equation*}
Furthermore, with
\begin{equation*}
a_k 
\triangleq \left< \hat{f} , e^{ik\pi\cdot} \right>_{L^2(-1,1)} 
= \dfrac{1}{2} \int_{-1}^1 f(\xi) e^{-i\tilde{k}\pi\xi} \diff\xi ,
\end{equation*}
and $\tilde{\overbrace{(-k-1)}^{}}=-\tilde{k}$, we obtain that
\begin{equation*}
f = \sum\limits_{k \in \mathbb{N}} a_k e^{i\tilde{k}\pi\cdot} + a_{-k-1} e^{-i\tilde{k}\pi\cdot} .
\end{equation*}

Let an arbitrary function $g \in L^2(0,1)$ be given and consider the functions $f_\mathrm{even},f_\mathrm{odd} \in L^2(-1,1)$ defined by
\begin{equation*}
f_\mathrm{even}(x) = 
\left\{\begin{split}
g(x) & \;\mathrm{if}\; x \geq 0 ; \\
g(-x) & \;\mathrm{if}\; x<0 .
\end{split}\right.
\;\;
f_\mathrm{odd}(x) = 
\left\{\begin{split}
g(x) & \;\mathrm{if}\; x \geq 0 ; \\
-g(-x) & \;\mathrm{if}\; x<0 .
\end{split}\right.
\end{equation*}

As $f_\mathrm{even}$ is an even function, we have,
\begin{equation*}
a_k = \dfrac{1}{2} \int_{-1}^1 f_\mathrm{even}(\xi) \cos(\tilde{k}\pi\xi) \diff\xi = a_{-k-1} ,
\end{equation*}
and thus
\begin{equation*}
f_\mathrm{even} = 2 \sum\limits_{k \in \mathbb{N}} a_k \cos(\tilde{k}\pi\cdot) .
\end{equation*}
As the above equality holds in $L^2(-1,1)$, it also holds in $L^2(0,1)$. Noting that $\left. f_\mathrm{even} \right|_{(0,1)} = g$, we conclude that $\{ \cos(\tilde{k}\pi\cdot) , \; k\in\mathbb{N}\}$ is maximal in $L^2(0,1)$.

Similarly, as $f_\mathrm{odd}$ is an odd function,
\begin{equation*}
a_k = - \dfrac{i}{2}  \int_{-1}^1 f_\mathrm{odd}(\xi) \sin(\tilde{k}\pi\xi) \diff\xi = - a_{-k-1} ,
\end{equation*}
and thus
\begin{equation*}
f = 2i \sum\limits_{k \in \mathbb{N}} a_k \sin(\tilde{k}\pi\cdot) .
\end{equation*}
Applying the same argument as above, we conclude that $\{ \sin(\tilde{k}\pi\cdot) , \; k\in\mathbb{N}\}$ is maximal in $L^2(0,1)$. \qed


\ifCLASSOPTIONcaptionsoff
  \newpage
\fi



\bibliographystyle{IEEEtranS}
\nocite{*}
\bibliography{IEEEabrv,mybibfile}

\begin{thebibliography}{10}
\providecommand{\url}[1]{#1}
\csname url@samestyle\endcsname
\providecommand{\newblock}{\relax}
\providecommand{\bibinfo}[2]{#2}
\providecommand{\BIBentrySTDinterwordspacing}{\spaceskip=0pt\relax}
\providecommand{\BIBentryALTinterwordstretchfactor}{4}
\providecommand{\BIBentryALTinterwordspacing}{\spaceskip=\fontdimen2\font plus
\BIBentryALTinterwordstretchfactor\fontdimen3\font minus
  \fontdimen4\font\relax}
\providecommand{\BIBforeignlanguage}[2]{{%
\expandafter\ifx\csname l@#1\endcsname\relax
\typeout{** WARNING: IEEEtranS.bst: No hyphenation pattern has been}%
\typeout{** loaded for the language `#1'. Using the pattern for}%
\typeout{** the default language instead.}%
\else
\language=\csname l@#1\endcsname
\fi
#2}}
\providecommand{\BIBdecl}{\relax}
\BIBdecl

\bibitem{argomedo2013strict}
F.~B. Argomedo, C.~Prieur, E.~Witrant, and S.~Br{\'e}mond, ``A strict control
  {L}yapunov function for a diffusion equation with time-varying distributed
  coefficients,'' \emph{IEEE Trans. Autom. Control}, vol.~58, no.~2, pp.
  290--303, Feb. 2013.

\bibitem{argomedo2012d}
F.~B. Argomedo, E.~Witrant, and C.~Prieur, ``{$D^1$}-{I}nput-to-state stability
  of a time-varying nonhomogeneous diffusive equation subject to boundary
  disturbances,'' in \emph{2012 American Control Conference}, 2012, pp.
  2978--2983.

\bibitem{bialy2016adaptive}
B.~J. Bialy, I.~Chakraborty, S.~C. Cekic, and W.~E. Dixon, ``Adaptive boundary
  control of store induced oscillations in a flexible aircraft wing,''
  \emph{Automatica}, vol.~70, pp. 230--238, 2016.

\bibitem{christensen2016introduction}
O.~Christensen \emph{et~al.}, \emph{{An Introduction to Frames and Riesz
  Bases}}.\hskip 1em plus 0.5em minus 0.4em\relax Springer, 2016.

\bibitem{Curtain2012}
R.~F. Curtain and H.~Zwart, \emph{An Introduction to Infinite-Dimensional
  Linear Systems Theory}.\hskip 1em plus 0.5em minus 0.4em\relax Springer
  Science \& Business Media, 2012, vol.~21.

\bibitem{dashkovskiy2011local}
S.~Dashkovskiy and A.~Mironchenko, ``Local {ISS} of reaction-diffusion
  systems,'' \emph{IFAC Proceedings Volumes}, vol.~44, no.~1, pp.
  11\,018--11\,023, 2011.

\bibitem{dashkovskiy2013input}
------, ``Input-to-state stability of infinite-dimensional control systems,''
  \emph{Math. Control Signals Syst.}, vol.~25, no.~1, pp. 1--35, 2013.

\bibitem{diagne2015state}
M.~Diagne and M.~Krstic, ``State-dependent input delay-compensated bang-bang
  control: Application to {3D} printing based on screw-extruder,'' in
  \emph{2015 American Control Conference}, 2015, pp. 5653--5658.

\bibitem{endo2017boundary}
T.~Endo, F.~Matsuno, and Y.~Jia, ``Boundary cooperative control by flexible
  {Timoshenko} arms,'' \emph{Automatica}, vol.~81, pp. 377--389, 2017.

\bibitem{guo2010spectral}
B.-Z. Guo, J.-M. Wang, and G.-D. Zhang, ``Spectral analysis of a wave equation
  with {K}elvin-{V}oigt damping,'' \emph{ZAMM Z. Angew. Math. Mech.}, vol.~90,
  no.~4, pp. 323--342, 2010.

\bibitem{guo2012spectrum}
B.-Z. Guo and G.-D. Zhang, ``On spectrum and {Riesz} basis property for
  one-dimensional wave equation with {Boltzmann} damping,'' \emph{ESAIM Control
  Optim. Calc. Var}, vol.~18, no.~3, pp. 889--913, 2012.

\bibitem{Hardy1952}
G.~H. Hardy, J.~E. Littlewood, and G.~P{\'o}lya, \emph{Inequalities}.\hskip 1em
  plus 0.5em minus 0.4em\relax Cambridge university press, 1952.

\bibitem{heigel2015thermo}
J.~Heigel, P.~Michaleris, and E.~Reutzel, ``Thermo-mechanical model development
  and validation of directed energy deposition additive manufacturing of
  {Ti-6Al-4V},'' \emph{Additive Manufacturing}, vol.~5, pp. 9--19, 2015.

\bibitem{henikl2016infinite}
J.~Henikl, W.~Kemmetm{\"u}ller, T.~Meurer, and A.~Kugi, ``Infinite-dimensional
  decentralized damping control of large-scale manipulators with hydraulic
  actuation,'' \emph{Automatica}, vol.~63, pp. 101--115, 2016.

\bibitem{jacob2018infinite}
B.~Jacob, R.~Nabiullin, J.~R. Partington, and F.~L. Schwenninger,
  ``Infinite-dimensional input-to-state stability and {Orlicz} spaces,''
  \emph{SIAM J. Control Optim.}, vol.~56, no.~2, pp. 868--889, 2018.

\bibitem{jacob2018continuity}
B.~Jacob, F.~L. Schwenninger, and H.~Zwart, ``On continuity of solutions for
  parabolic control systems and input-to-state stability,'' \emph{J.
  Differential Equations}, 2018.

\bibitem{karafyllis2018boundary}
I.~Karafyllis, M.~Kontorinaki, and M.~Krstic, ``Boundary-to-displacement
  asymptotic gains for wave systems with {Kelvin-Voigt} damping,'' \emph{arXiv
  preprint arXiv:1807.06549}, 2018.

\bibitem{karafyllis2016input}
I.~Karafyllis and M.~Krstic, ``Input-to state stability with respect to
  boundary disturbances for the {1-D} heat equation,'' in \emph{IEEE 55th
  Conference on Decision and Control}, 2016, pp. 2247--2252.

\bibitem{karafyllis2017iss}
------, ``{ISS} in different norms for {1-D} parabolic {PDEs} with boundary
  disturbances,'' \emph{SIAM J. Control Optim.}, vol.~55, no.~3, pp.
  1716--1751, 2017.

\bibitem{karafyllis2018small}
------, ``Small-gain stability analysis of certain hyperbolic--parabolic {PDE}
  loops,'' \emph{Systems \& Control Letters}, vol. 118, pp. 52--61, 2018.

\bibitem{karafyllis2019preview}
------, \emph{Input-to-State Stability for {PDEs}}.\hskip 1em plus 0.5em minus
  0.4em\relax Springer, 2019.

\bibitem{karafyllis2016iss}
------, ``{ISS} with respect to boundary disturbances for {1-D} parabolic
  {PDEs},'' \emph{IEEE Trans. Autom. Control}, vol.~61, no.~12, pp. 3712--3724,
  Dec. 2016.

\bibitem{kato2013perturbation}
T.~Kato, \emph{Perturbation Theory for Linear Operators}.\hskip 1em plus 0.5em
  minus 0.4em\relax Springer Science \& Business Media, 2013, vol. 132.

\bibitem{korner1989fourier}
T.~W. K{\"o}rner, \emph{Fourier Analysis}.\hskip 1em plus 0.5em minus
  0.4em\relax Cambridge university press, 1989.

\bibitem{lhachemi2018boundaryAutomatica}
H.~Lhachemi, D.~Saussi{\'e}, and G.~Zhu, ``Boundary feedback stabilization of a
  flexible wing model under unsteady aerodynamic loads,'' \emph{Automatica},
  vol.~97, pp. 73--81, 2018.

\bibitem{lhachemi2018boundary}
------, ``Boundary control of a nonhomogeneous flexible wing with bounded input
  disturbances,'' \emph{IEEE Trans. Autom. Control}, vol.~64, no.~2, pp.
  854--861, Feb. 2019.

\bibitem{liu2006exponential}
K.~Liu and B.~Rao, ``Exponential stability for the wave equations with local
  {K}elvin–{V}oigt damping,'' \emph{Zeitschrift für angewandte Mathematik
  und Physik}, vol.~57, no.~3, pp. 419--432, May 2006.

\bibitem{luo2012stability}
Z.-H. Luo, B.-Z. Guo, and {\"O}.~Morg{\"u}l, \emph{Stability and Stabilization
  of Infinite Dimensional Systems with Applications}.\hskip 1em plus 0.5em
  minus 0.4em\relax Springer Science \& Business Media, 2012.

\bibitem{mazenc2011strict}
F.~Mazenc and C.~Prieur, ``Strict {L}yapunov functions for semilinear parabolic
  partial differential equations,'' \emph{Mathematical Control and Related
  Fields}, vol.~1, no.~2, pp. 231--250, 2011.

\bibitem{mironchenko2016local}
A.~Mironchenko, ``Local input-to-state stability: Characterizations and
  counterexamples,'' \emph{Systems \& Control Letters}, vol.~87, pp. 23--28,
  Jan. 2016.

\bibitem{mironchenko2015construction}
A.~Mironchenko and H.~Ito, ``Construction of {L}yapunov functions for
  interconnected parabolic systems: an {iISS} approach,'' \emph{SIAM J. Control
  Optim.}, vol.~53, no.~6, pp. 3364--3382, 2015.

\bibitem{mironchenko2017monotonicity}
A.~Mironchenko, I.~Karafyllis, and M.~Krstic, ``Monotonicity methods for
  input-to-state stability of nonlinear parabolic {PDEs} with boundary
  disturbances,'' \emph{arXiv preprint arXiv:1706.07224}, 2017.

\bibitem{mironchenko2016restatements}
A.~Mironchenko and F.~Wirth, ``Restatements of input-to-state stability in
  infinite dimensions: what goes wrong,'' in \emph{Proc. of 22th International
  Symposium on Mathematical Theory of Systems and Networks}, 2016.

\bibitem{mironchenko2017characterizations}
------, ``Characterizations of input-to-state stability for
  infinite-dimensional systems,'' \emph{IEEE Trans. Autom. Control}, vol.~63,
  no.~6, pp. 1692--1707, June 2018.

\bibitem{Pazy2012}
A.~Pazy, \emph{Semigroups of Linear Operators and Applications to Partial
  Differential Equations}.\hskip 1em plus 0.5em minus 0.4em\relax Springer
  Science \& Business Media, 2012, vol.~44.

\bibitem{prieur2012iss}
C.~Prieur and F.~Mazenc, ``{ISS}-{L}yapunov functions for time-varying
  hyperbolic systems of balance laws,'' \emph{Math. Control Signals Syst.},
  vol.~24, no. 1-2, pp. 111--134, 2012.

\bibitem{sontag1989smooth}
E.~D. Sontag, ``Smooth stabilization implies coprime factorization,''
  \emph{IEEE Trans. Autom. Control}, vol.~34, no.~4, pp. 435--443, Apr. 1989.

\bibitem{tanwani2017disturbance}
A.~Tanwani, C.~Prieur, and S.~Tarbouriech, ``Disturbance-to-state stabilization
  and quantized control for linear hyperbolic systems,'' \emph{arXiv preprint
  arXiv:1703.00302}, 2017.

\bibitem{xu2008spectrum}
G.~Q. Xu and N.~E. Mastorakis, ``Spectrum of an operator arising elastic system
  with local {K-V} damping,'' \emph{ZAMM Z. Angew. Math. Mech.}, vol.~88,
  no.~6, pp. 483--496, 2008.

\bibitem{zheng2018iss}
J.~Zheng, H.~Lhachemi, G.~Zhu, and D.~Saussi{\'e}, ``{ISS} with respect to
  boundary and in-domain disturbances for a coupled beam-string system,''
  \emph{Math. Control Signals Systems}, vol.~30, no.~4, p.~21, 2018.

\bibitem{zheng2017input}
J.~Zheng and G.~Zhu, ``Input-to-state stability with respect to boundary
  disturbances for a class of semi-linear parabolic equations,''
  \emph{Automatica}, vol.~97, pp. 271--277, 2018.

\bibitem{zheng2017giorgi}
------, ``{A De Giorgi} iteration-based approach for the establishment of {ISS}
  properties for {B}urgers' equation with boundary and in-domain
  disturbances,'' \emph{IEEE Trans. Autom. Control}, in press, 2018.

\end{thebibliography}

\end{document}